\documentclass[11pt,reqno]{amsart}
\usepackage{amssymb}
\usepackage{amsmath}
\usepackage{amscd}
\usepackage[hidelinks]{hyperref}
\begin{document}
\setlength{\oddsidemargin}{0cm} \setlength{\evensidemargin}{0cm}

\theoremstyle{plain}
\newtheorem{theorem}{Theorem}[section]
\newtheorem{proposition}[theorem]{Proposition}
\newtheorem{lemma}[theorem]{Lemma}
\newtheorem{corollary}[theorem]{Corollary}
\newtheorem{conj}[theorem]{Conjecture}

\theoremstyle{definition}
\newtheorem{definition}[theorem]{Definition}
\newtheorem{exam}[theorem]{Example}
\newtheorem{remark}[theorem]{Remark}

\numberwithin{equation}{section}

\title[The decomposition of a Lie group and the uniqueness]
{The decomposition of a Lie group with a left invariant pseudo-Riemannian metric and the uniqueness}

\author{Zhiqi Chen}
\address{School of Mathematical Sciences and LPMC \\ Nankai University \\ Tianjin 300071, P.R. China} \email{chenzhiqi@nankai.edu.cn}

\author{Ke Liang}
\address{School of Mathematical Sciences and LPMC \\ Nankai University \\ Tianjin 300071, P.R. China}\email{liangke@nankai.edu.cn}

\author{Mingming Ren}
\address{Beijing Computational Science Research Center \\ No. 3, Heqing Road \\ Haidian District \\ Beijing 100084, P.R. China}\email{rmingming@gmail.com}



\begin{abstract}
In this paper, we discuss the decomposition of a Lie group with a
left invariant pseudo-Riemannian metric and the uniqueness. The
decomposition is defined by non-degenerate strong ideals of the
Levi-Civita connection. Its factors are totally geodesic normal
subgroups on the simply connected covering group. We explain its
relation with the de Rham--Wu decomposition and with the usual
decompositions of Lie algebras by explicit examples, including a
non-flat Riemannian product which is indecomposable in the present
sense. We also describe the possible mixed metric terms between
strong factors. As applications, we obtain the factorwise reduction
of the Einstein equation and of the geodesic equation. In particular,
the strong factors are unique up to their order when the Ricci tensor
is non-degenerate, and geodesic completeness reduces to completeness
of the factors even when the strong decomposition is not orthogonal.
\end{abstract}

\maketitle

\setcounter{section}{-1}
\section{Introduction}
Lie groups play an important role in modern geometry. A left
invariant pseudo-Riemannian metric on a Lie group gives both an
algebraic structure and a geometric structure. It is natural to ask
which decompositions preserve these two structures at the same time.
The de Rham theorem and its pseudo-Riemannian extension due to Wu
describe local orthogonal products in terms of non-degenerate
parallel distributions \cite{de,Wu1,Wu2}. Their global versions
require additional hypotheses, in particular simple connectivity
and geodesic completeness. There are many studies on these
decompositions and on pseudo-Riemannian holonomy
\cite{Ber,Ch,GT,Ma}.
\medskip

On a Lie group, a parallel distribution need not give a normal
subgroup. Even when it is left invariant, the corresponding
subspace of the Lie algebra need not be a Lie ideal. Thus the
decomposition of the metric manifold and the decomposition compatible
with group multiplication are different questions. This distinction
already occurs for Riemannian metrics. We give a flat example on the
universal covering group of the Euclidean motion group and a
non-flat example isometric to $\mathbb H^2\times\mathbb R^2$.
Both are indecomposable in the sense used in this paper.
\medskip

Let $G$ be a connected Lie group, let $\mathfrak g$ be its Lie algebra,
and let $\langle,\rangle$ be a left invariant pseudo-Riemannian
metric. We use the Levi-Civita product $\nabla_XY$ on $\mathfrak g$.
A subspace $\mathfrak h$ is called a strong ideal if
$$\nabla_{\mathfrak g}\mathfrak h\subseteq\mathfrak h,
\qquad \nabla_{\mathfrak h}\mathfrak g\subseteq\mathfrak h.$$
We discuss decompositions into indecomposable, non-degenerate strong
ideals. Such a decomposition gives a direct product of Lie groups on
the simply connected cover, and its mixed covariant derivatives
vanish. The metric terms between distinct factors, however, need
not vanish. We describe these terms by bilinear forms on the
quotients $\mathfrak g_i/\nabla_{\mathfrak g_i}\mathfrak g_i$.
This gives a precise measure of the possible failure of
orthogonality.
\medskip

For a bi-invariant metric, strong ideals are exactly Lie ideals,
and the orthogonal complement of a non-degenerate ideal is again an
ideal. Decomposition results in this setting were considered in
\cite{ZZ}. For a general left invariant metric, the orthogonal
complement is parallel but need not be a strong ideal. Related
algebraic questions were considered in \cite{CLZ,CHB}.
\medskip

The right annihilator $Ann_R(\mathfrak g)$ consists of the left
invariant parallel vector fields. Its relation with
$\nabla_{\mathfrak g}\mathfrak g$ is central to the uniqueness
problem. We prove uniqueness up to factorwise strong isometries when
$Ann_R(\mathfrak g)$ is totally isotropic. When
$Ann_R(\mathfrak g)=Ann(\mathfrak g)$, an orthogonal decomposition
exists and is unique up to a strong isometry, with a maximal
non-degenerate annihilator summand grouped together. If
$Ann_R(\mathfrak g)=0$, the factors themselves are unique up to
their order. The last condition follows from non-degeneracy of the
Ricci tensor, but its converse does not hold.
\medskip

The Einstein equation therefore reduces to the indecomposable
strong factors when the Einstein constant is non-zero. This is a
reduction for a given metric; it does not assert that every Einstein
metric on an algebraic direct product must split along prescribed
Lie ideals. Other decomposition and existence results for Einstein
metrics can be found in \cite{Bo,Ja1,Ja2,Ni1,Ni2}. General references
include \cite{Be,An,Ber1,Ber2,LW}, and for homogeneous Einstein
metrics see \cite{DZ1,GK,He,La1,La2}. A second application is the
reduction of the geodesic equation and of geodesic completeness,
which also applies to non-orthogonal strong decompositions.
\medskip

The paper is organized as follows. Section~\ref{sect1} gives the
definitions and some facts about the Levi-Civita connection.
Section~\ref{sect2} gives the local and global decompositions,
an algebraic criterion for finding the factors, and examples
comparing them with other decompositions.
Section~\ref{sect3} discusses uniqueness according to the structure
of the annihilator. Section~\ref{sect4} gives the applications to
Einstein metrics and geodesics. The final section discusses the
meaning and the limits of the decomposition results.

\section{Preliminary}\label{sect1}
Throughout this paper, Lie algebras are real and finite-dimensional,
and Lie groups are connected. Let $G$ be a Lie group, let
$\mathfrak g$ consist of its left invariant vector fields, and let
$\langle,\rangle$ be a left invariant pseudo-Riemannian metric.
The unique torsion-free metric connection, namely the Levi-Civita
connection, is determined by
\begin{equation}\label{nabla1}
2\langle\nabla_XY,Z\rangle
=\langle[X,Y],Z\rangle-\langle[Y,Z],X\rangle
+\langle[Z,X],Y\rangle.
\end{equation}
Equivalently, for any $X,Y,Z\in\mathfrak g$,
\begin{align}
\nabla_XY-\nabla_YX&=[X,Y],\label{torsion}\\
\langle\nabla_XY,Z\rangle+\langle Y,\nabla_XZ\rangle&=0.
\label{metriccompat}
\end{align}

A subspace $\mathfrak h$ is called {\bf totally isotropic} if
$\langle X,Y\rangle=0$ for all $X,Y\in\mathfrak h$. We also use
the shorter term {\bf isotropic} with this meaning. It is called
{\bf non-degenerate} if the restriction of the metric to it is
non-degenerate. Let $\nabla_{\mathfrak h_1}\mathfrak h_2$ denote the
linear span of all $\nabla_XY$, where $X\in\mathfrak h_1$ and
$Y\in\mathfrak h_2$. A subspace $\mathfrak h$ is a {\bf strong
ideal} if both $\nabla_{\mathfrak g}\mathfrak h$ and
$\nabla_{\mathfrak h}\mathfrak g$ are contained in $\mathfrak h$.
It is a Lie ideal by (\ref{torsion}).
\medskip

A linear map $\pi:\mathfrak g\longrightarrow\mathfrak g'$ is a
{\bf strong homomorphism} if
$$\pi(\nabla_XY)=\nabla'_{\pi(X)}\pi(Y).$$
It is a {\bf strong isomorphism} if it is also bijective, and a
{\bf strong automorphism} if $\mathfrak g=\mathfrak g'$. A strong
isomorphism is a {\bf strong isometry} if
$$\langle\pi(X),\pi(Y)\rangle'=\langle X,Y\rangle.$$
Every strong homomorphism is a Lie algebra homomorphism. Conversely,
an isometric Lie algebra isomorphism is strong by (\ref{nabla1}).
For Lie groups, a strong homomorphism is a Lie group homomorphism
whose differential is strong; a strong isomorphism is required to
be a Lie group isomorphism as well. A bijective differential alone
does not imply a group isomorphism. The corresponding terminology
for strong isometries is used in the same way.
\medskip

Write
\begin{align*}
\mathfrak h^\perp&=\{X\in\mathfrak g:\langle X,Y\rangle=0
\text{ for all }Y\in\mathfrak h\},\\
Ann_R(\mathfrak g)&=\{X\in\mathfrak g:\nabla_YX=0
\text{ for all }Y\in\mathfrak g\},\\
Ann(\mathfrak g)&=\{X\in\mathfrak g:\nabla_YX=\nabla_XY=0
\text{ for all }Y\in\mathfrak g\}.
\end{align*}
The fields in $Ann_R(\mathfrak g)$ are parallel and hence Killing.
Indeed, vanishing of their covariant derivatives in a left invariant
frame implies vanishing in every direction. Also
$Ann(\mathfrak g)$ is central, by (\ref{torsion}).

\begin{proposition}\label{prop2.1}
$Ann_R(\mathfrak g)=(\nabla_{\mathfrak g}\mathfrak g)^\perp$.
\end{proposition}
\begin{proof}
For $X\in\mathfrak g$, equation (\ref{metriccompat}) gives
$$\langle\nabla_YX,Z\rangle=-\langle X,\nabla_YZ\rangle.$$
The left side vanishes for all $Y,Z$ if and only if the right side
does. Non-degeneracy of the metric proves the assertion.
\end{proof}

\begin{proposition}\label{prop1}
If $\mathfrak h$ is a strong ideal, then
$\nabla_{\mathfrak g}\mathfrak h^\perp\subseteq\mathfrak h^\perp$
and $\nabla_{\mathfrak h}\mathfrak h^\perp=0$.
\end{proposition}
\begin{proof}
For $X\in\mathfrak g$, $Y\in\mathfrak h$ and
$Z\in\mathfrak h^\perp$, we have
$$\langle Y,\nabla_XZ\rangle=-\langle\nabla_XY,Z\rangle=0,
\qquad
\langle X,\nabla_YZ\rangle=-\langle\nabla_YX,Z\rangle=0.$$
These two identities give the assertions.
\end{proof}

\section{The decomposition and geometry}\label{sect2}
\subsection{Strong factors and the metric}
\begin{definition}
The metric Lie algebra $\mathfrak g$ is called {\bf decomposable}
if $\mathfrak g=\mathfrak h_1\oplus\mathfrak h_2$ for two
non-zero, non-degenerate strong ideals. Otherwise it is called
{\bf indecomposable}. The decomposition is {\bf orthogonal} if
$\langle\mathfrak h_1,\mathfrak h_2\rangle=0$. The same terminology
is used for the corresponding metric Lie group.
\end{definition}

\begin{remark}\label{rem1}
A non-zero proper non-degenerate strong ideal $\mathfrak h$ gives
an orthogonal vector space decomposition
$\mathfrak g=\mathfrak h\oplus\mathfrak h^\perp$ into parallel
subspaces. But $\mathfrak h^\perp$ need not be a strong ideal.
Consequently the existence of $\mathfrak h$ does not imply
decomposability in the present sense. Example~\ref{ex:e2} below
gives an explicit instance of this distinction.
\end{remark}

\begin{proposition}\label{prop2}
There exists a decomposition
$\mathfrak g=\mathfrak g_1\oplus\cdots\oplus\mathfrak g_l$
into indecomposable, non-degenerate strong ideals.
\end{proposition}
\begin{proof}
If a factor is decomposable, split it further. Mixed covariant
derivatives between strong summands vanish, so a strong ideal of a
summand remains a strong ideal of $\mathfrak g$. The procedure
terminates because dimension decreases at each splitting.
\end{proof}

\begin{proposition}\label{prop:mixed}
Suppose $\mathfrak g=\bigoplus_i\mathfrak g_i$ is a direct sum of
strong ideals. Then, for $i\ne j$,
\begin{equation}\label{eq:mixed}
\nabla_{\mathfrak g_i}\mathfrak g_j=0,
\qquad [\mathfrak g_i,\mathfrak g_j]=0,
\qquad
\langle\nabla_{\mathfrak g_i}\mathfrak g_i,\mathfrak g_j\rangle=0.
\end{equation}
If the restrictions of the metric are non-degenerate, the induced
Levi-Civita connections are the restrictions of $\nabla$.
\end{proposition}
\begin{proof}
The first mixed derivative belongs to both ideals and hence is
zero. The bracket assertion follows from (\ref{torsion}). For
$X,Y\in\mathfrak g_i$ and $Z\in\mathfrak g_j$,
$\langle\nabla_XY,Z\rangle=-\langle Y,\nabla_XZ\rangle=0$.
Finally the restricted connection is torsion-free and metric.
\end{proof}

\begin{proposition}\label{prop:coupling}
Let $(\mathfrak g_i,\langle,\rangle_i)$, $i=1,2$, be metric Lie
algebras with Levi-Civita connections $\nabla^i$, and put
$D_i=\nabla^i_{\mathfrak g_i}\mathfrak g_i$. Let $B$ be a bilinear
form on $\mathfrak g_1\times\mathfrak g_2$, and suppose that
\begin{align*}
\langle X_1+X_2,Y_1+Y_2\rangle_B
={}&\langle X_1,Y_1\rangle_1+\langle X_2,Y_2\rangle_2\\
&+B(X_1,Y_2)+B(Y_1,X_2)
\end{align*}
is non-degenerate on the direct sum Lie algebra. Then both summands
are strong ideals if and only if
\begin{equation}\label{eq:coupling}
B(D_1,\mathfrak g_2)=0,
\qquad B(\mathfrak g_1,D_2)=0.
\end{equation}
In this case the connection is $\nabla^1\oplus\nabla^2$.
\end{proposition}
\begin{proof}
Necessity follows from Proposition~\ref{prop:mixed}. Conversely,
the direct sum connection is torsion-free. It preserves the
restrictions of the metric. Its metric compatibility for arguments
in different summands is exactly (\ref{eq:coupling}); the other
mixed terms vanish. It is therefore the Levi-Civita connection of
$\langle,\rangle_B$.
\end{proof}

Thus the possible mixed forms are precisely the bilinear forms on
$(\mathfrak g_1/D_1)\times(\mathfrak g_2/D_2)$ for which the full
metric is non-degenerate. Before imposing non-degeneracy, this
space has dimension
$$\dim Ann_R(\mathfrak g_1)\,\dim Ann_R(\mathfrak g_2).$$
In particular, varying such a mixed form does not change the
connection or its curvature tensor of type $(1,3)$.

\begin{proposition}\label{prop:projection}
For $X\in\mathfrak g$, set $L_XY=\nabla_XY$ and $T_XY=\nabla_YX$.
A linear map $P$ is the projection associated with a direct sum of
strong ideals if and only if
\begin{equation}\label{eq:projector}
P^2=P,\qquad PL_X=L_XP,\qquad PT_X=T_XP
\quad\text{for all }X\in\mathfrak g.
\end{equation}
It gives an orthogonal decomposition if and only if it is also
self-adjoint. To give a decomposition in the sense of the definition,
both its image and its kernel must be non-zero and non-degenerate.
\end{proposition}
\begin{proof}
Commutation with $L_X$ and $T_X$ is equivalent to invariance of the
image and kernel under all left and right multiplications. For an
idempotent, self-adjointness is equivalent to orthogonality of the
image and kernel; non-degeneracy of each then follows from that of
$\mathfrak g$.
\end{proof}

In a basis, the commutation equations in (\ref{eq:projector}) are
linear equations for the entries of $P$. The additional equation
$P^2=P$ is quadratic. This gives a finite algebraic procedure for
detecting factors; it is not merely a linear test. In indefinite
signature one must also check the non-degeneracy conditions.

\subsection{Global decomposition and parallel distributions}
\begin{proposition}[\cite{Kn1}]\label{prop3}
Let $G$ be simply connected, and let
$\mathfrak g=\mathfrak h_1\oplus\mathfrak h_2$ be a direct sum of
Lie ideals. The corresponding connected analytic subgroups $H_i$
are closed, normal and simply connected, and multiplication gives
a Lie group isomorphism $H_1\times H_2\longrightarrow G$.
\end{proposition}

\begin{proposition}\label{prop2-1}
If $G$ is simply connected, a strong decomposition gives a Lie
group decomposition $G=G_1\times\cdots\times G_l$ into closed,
normal, simply connected analytic subgroups. Every $G_i$ is totally
geodesic. Multiplication identifies the Levi-Civita connection with
the direct sum connection; it is an isometry from the product of
the restricted metrics if and only if the decomposition is orthogonal.
\end{proposition}
\begin{proof}
Use Propositions~\ref{prop3} and \ref{prop:mixed}. In particular,
$\nabla_XY\in\mathfrak g_i$ for $X,Y\in\mathfrak g_i$, so the
second fundamental form of $G_i$ vanishes. The metric assertion
can be checked at the identity by left invariance.
\end{proof}

\begin{proposition}\label{prop2-2}
For a connected $G$, the preceding decomposition holds on its
simply connected covering group $\widetilde G$, equipped with the
lifted metric.
\end{proposition}

\begin{remark}\label{rem:descent}
Write $G=\widetilde G/\Gamma$, where $\Gamma$ is discrete and
central, and identify $\widetilde G=\prod_i\widetilde G_i$.
The chosen product descends to a product of the corresponding
quotient factors if and only if
$$\Gamma=\prod_i(\Gamma\cap\widetilde G_i).$$
There is no such conclusion for an arbitrary $\Gamma$. For example,
in the flat torus $\mathbb R^2/\mathbb Z^2$, an orthogonal splitting
along lines of irrational slope gives a strong splitting of its
Lie algebra, but the corresponding analytic subgroups are not
closed. The covering-group statement avoids this obstruction.
\end{remark}

\begin{remark}\label{rem:derham}
The left invariant distribution associated with a strong ideal is
parallel. For a non-degenerate strong ideal its orthogonal
complement is also parallel, by Proposition~\ref{prop1}. Hence it
gives a local pseudo-Riemannian product in the sense of de
Rham--Wu. The additional requirement in this paper concerns
compatibility with Lie ideals in both factors. A global product of
simply connected Lie groups in Proposition~\ref{prop2-1} does not
require geodesic completeness. A global de Rham--Wu theorem for
an arbitrary pseudo-Riemannian manifold does require further
hypotheses \cite{Wu1,Wu2}. Also, a null parallel subspace alone is
not a non-degenerate factor and does not imply a metric product.
\end{remark}

\subsection{Examples and comparison with other decompositions}
In the examples below, brackets and scalar products not specified
are zero, apart from those forced by skew-symmetry or symmetry.

\begin{exam}\label{ex:e2}
Let $\mathfrak e$ have orthonormal basis $\{a,u,v\}$ and brackets
$$[a,u]=v,\qquad [a,v]=-u.$$
This is the Lie algebra of the orientation preserving Euclidean
motion group. By (\ref{nabla1}), the only non-zero products are
$$\nabla_a u=v,\qquad \nabla_a v=-u.$$
Consequently the metric is flat, $Ann_R(\mathfrak e)=\mathbb Ra$
and $Ann(\mathfrak e)=0$. On the simply connected group it is the
complete Euclidean metric in coordinates $(t,x,y)$: left translation
rotates the $(x,y)$ coordinates and preserves $dt^2+dx^2+dy^2$.
This is also a basic example in Milnor's description of flat Lie
groups \cite{Mi}.
\medskip

Nevertheless $\mathfrak e$ is indecomposable. A non-trivial direct
sum of ideals in dimension three would have a one-dimensional
summand. Such a summand would be central, whereas the center of
$\mathfrak e$ is zero. The plane $\mathfrak u=\operatorname{span}
\{u,v\}$ is a non-degenerate strong ideal, but
$\mathfrak u^\perp=\mathbb Ra$ is not a Lie ideal. Thus there is a
non-trivial de Rham metric product $\mathbb R\times\mathbb R^2$,
although there is no non-trivial strong decomposition. The usual
de Rham convention groups all these Euclidean directions into one
flat factor; the splitting of that flat factor into lines is not
canonical.
\end{exam}

\begin{exam}\label{ex:twisted}
There is a non-flat example with $Ann_R(\mathfrak g)=0$. Let
$\{a,u,v,w\}$ be orthonormal, with
$$[a,u]=u,\qquad [a,v]=w,\qquad [a,w]=-v.$$
The only non-zero products are
\begin{equation}\label{eq:twistedconnection}
\nabla_u a=-u,\quad \nabla_u u=a,\quad
\nabla_a v=w,\quad \nabla_a w=-v.
\end{equation}
The simply connected group has multiplication
$$ (t,x,z)(s,y,\zeta)
=(t+s,x+e^t y,z+R_t\zeta),$$
where $z,\zeta\in\mathbb R^2$ and $R_t$ is rotation through angle
$t$. Its left invariant metric is
$$dt^2+e^{-2t}dx^2+dz_1^2+dz_2^2.$$
Thus the metric manifold is globally $\mathbb H^2\times\mathbb R^2$,
with the hyperbolic factor of curvature $-1$.
\medskip

We show that its strong decomposition is trivial. Any projection
$P$ in Proposition~\ref{prop:projection} preserves
$\operatorname{span}\{a,u\}$ and $\operatorname{span}\{v,w\}$ by
commuting with $L_u$ and $L_a$. On the first plane write
$Pa=\alpha a+\beta u$. Commutation with $L_u$ determines $Pu$.
The identity $P\nabla_a a=\nabla_{Pa}a$ gives $\beta=0$, so
$Pa=\alpha a$ and $Pu=\alpha u$. Next
$P\nabla_a v=\nabla_{Pa}v$ and
$P\nabla_a w=\nabla_{Pa}w$ give $Pv=\alpha v$ and $Pw=\alpha w$.
Hence $P=\alpha I$, and idempotence gives $P=0$ or $I$.
\medskip

Both planes are parallel, but the first one is not a Lie ideal.
Moreover (\ref{eq:twistedconnection}) gives
$$Ann_R(\mathfrak g)=0,\qquad
ric=\operatorname{diag}(-1,-1,0,0).$$
In particular, vanishing of the right annihilator does not imply
non-degeneracy of the Ricci tensor. The Euclidean factor has
parallel vector fields, but none of its non-zero parallel fields
is left invariant under this transitive group action.
\end{exam}

\begin{exam}\label{ex:compact}
Let $\mathfrak k=\mathfrak{su}(2)$ with an invariant positive
definite form $Q$, and let $0<|\rho|<1$. On
$\mathfrak g=\mathfrak k\oplus\mathfrak k$ put
\begin{align*}
\langle(x_1,x_2),(y_1,y_2)\rangle
={}&Q(x_1,y_1)+Q(x_2,y_2)\\
&+\rho\{Q(x_1,y_2)+Q(x_2,y_1)\}.
\end{align*}
Both simple ideals are non-degenerate and both corresponding
normal subgroups are totally geodesic. Indeed, on either ideal
$\nabla_XY=\frac12[X,Y]$, as follows by pairing with vectors in
each summand in (\ref{nabla1}). However,
$$2\langle\nabla_{(x,0)}(0,y),(z,0)\rangle
=\rho Q([z,x],y),$$
which is not identically zero. Thus neither ideal is strong. The
only possible non-trivial direct decomposition of this semisimple
Lie algebra consists of its two simple ideals. It follows that the
metric Lie algebra is indecomposable in the present sense.
\medskip

This example distinguishes strong decomposition from the ordinary
direct decomposition into Lie ideals. It also shows that totally
geodesic normal factors need not be orthogonal, and that total
geodesy of the factors alone does not give a strong decomposition.
\end{exam}

\begin{exam}\label{ex:iwasawa}
For $\mathfrak{sl}(2,\mathbb R)$ write
$$H=\begin{pmatrix}1&0\\0&-1\end{pmatrix},\quad
E=\begin{pmatrix}0&1\\0&0\end{pmatrix},\quad
F=\begin{pmatrix}0&0\\1&0\end{pmatrix}.$$
Its Iwasawa vector space decomposition is
$$\mathfrak g=\mathbb R(E-F)\oplus\mathbb RH\oplus\mathbb RE
=\mathfrak k\oplus\mathfrak a\oplus\mathfrak n.$$
But $[H,E]=2E$ and $[E-F,E]=H$, so these are not three commuting
ideals. The Cartan decomposition
$\mathfrak g=\mathfrak k\oplus\operatorname{span}\{H,E+F\}$ is
also not a decomposition into ideals. Since $\mathfrak g$ is
simple, it is strongly indecomposable for every non-degenerate
metric. The group decomposition $KAN$ and the Cartan decomposition
therefore address a different question \cite{Kn1}.
\medskip

Similarly, a Levi decomposition
$\mathfrak g=\mathfrak s\ltimes\mathfrak r$ separates a semisimple
subalgebra and the solvable radical. It need not be a direct sum of
ideals or respect the metric. In Example~\ref{ex:e2}, the radical
is the whole algebra, while the non-trivial splitting
$\mathbb Ra\ltimes\operatorname{span}\{u,v\}$ is a semidirect
product and is not a strong direct decomposition.
\end{exam}

\subsection{A qualification concerning minimal subgroups}
\begin{proposition}\label{prop:minimal}
Let a Riemannian metric Lie algebra be an algebraic direct sum of
ideals $\mathfrak g=\mathfrak h\oplus\mathfrak k$, without an
orthogonality assumption. Let $H$ be the connected subgroup with
Lie algebra $\mathfrak h$. For an orthonormal basis $e_1,\ldots,e_r$
of $\mathfrak h$, its mean curvature vector $\mathcal H$ satisfies
\begin{equation}\label{eq:mean}
r\langle\mathcal H,Z\rangle
=\operatorname{tr}(\operatorname{ad}Z|_{\mathfrak h}),
\qquad Z\in\mathfrak h^\perp.
\end{equation}
In particular $H$ is minimal if $\mathfrak h$ is unimodular.
\end{proposition}
\begin{proof}
The Koszul formula gives
$\langle\nabla_{e_i}e_i,Z\rangle=\langle[Z,e_i],e_i\rangle$.
Since $\mathfrak h$ is an ideal, summing is the indicated trace.
Write $Z=Z_{\mathfrak h}+Z_{\mathfrak k}$. The ideals commute, so
this trace is $\operatorname{tr}(\operatorname{ad}Z_{\mathfrak h}
|_{\mathfrak h})$, which vanishes if $\mathfrak h$ is unimodular.
\end{proof}

\begin{exam}\label{ex:nonminimal}
Let $\mathfrak g=\operatorname{span}\{a,b\}\oplus\mathbb Rz$, with
$[a,b]=b$, and give it a metric for which $a,b,z$ have unit length
and $\langle a,z\rangle=\rho$, where $0<|\rho|<1$.
For $\mathfrak h=\operatorname{span}\{a,b\}$ the normal vector
$Z=z-\rho a$ satisfies
$$2\langle\mathcal H,Z\rangle
=\operatorname{tr}(\operatorname{ad}Z|_{\mathfrak h})=-\rho.$$
Thus this normal direct factor is not minimal. Example~\ref{ex:compact}
shows a separate issue: even for semisimple groups, total geodesy
of the two factors does not imply their orthogonality.
\end{exam}

\section{The uniqueness of the decomposition}\label{sect3}
We first discuss the decomposition of the metric Lie algebra. On
the simply connected group the corresponding statements follow by
integration. For other connected groups one must also take account
of the covering kernel, as in Remark~\ref{rem:descent}.

\subsection[The non-isotropic right annihilator]{The case when \texorpdfstring{$Ann_R(\mathfrak g)$}{the right annihilator} is not isotropic}
\label{subsect2}
If $Ann_R(\mathfrak g)=\mathfrak g$, the product is zero. There is
an orthogonal decomposition into one-dimensional non-degenerate
strong ideals, but the particular lines need not be unique.
\medskip

Suppose $Ann_R(\mathfrak g)$ is not isotropic. Choose a maximal
non-degenerate subspace $\mathfrak h_1$ of $Ann_R(\mathfrak g)$ and
put $\mathfrak g_1=\mathfrak h_1^\perp$. For $X\in\mathfrak g$,
$Y\in\mathfrak h_1$ and $Z\in\mathfrak g_1$,
$$\langle\nabla_XZ,Y\rangle=-\langle Z,\nabla_XY\rangle=0.$$
Thus $\nabla_{\mathfrak g}\mathfrak g_1\subseteq\mathfrak g_1$.
Also $\nabla_{\mathfrak g_1}\mathfrak h_1=0$, so $\mathfrak g_1$
is a strong ideal. The subspace $\mathfrak h_1$ is an abelian
subalgebra, since its right products vanish. We obtain
$$\mathfrak g=\mathfrak h_1\ltimes\mathfrak g_1.$$
It is not asserted that $\mathfrak h_1$ is a strong ideal.
Repeating this construction inside $\mathfrak g_1$ gives a finite
sequence ending in an algebra with isotropic right annihilator.
The final algebra is allowed to be zero.

\begin{theorem}\label{thm1}
Suppose $G$ is simply connected and $Ann_R(\mathfrak g)$ is not
isotropic. There is a sequence of closed, connected, simply
connected analytic subgroups
$$G=G_0\supset G_1\supset\cdots\supset G_n$$
such that each $G_i$ is normal and totally geodesic in $G_{i-1}$,
each quotient $G_{i-1}/G_i$ is an abelian vector group, and each
$\mathfrak g_i$ is a non-degenerate strong ideal of
$\mathfrak g_{i-1}$. Moreover $Ann_R(\mathfrak g_n)$ is isotropic.
For a connected group that is not simply connected, these
conclusions hold on its simply connected covering group.
\end{theorem}
\begin{proof}
Integrate each semidirect sum constructed above to the simply
connected semidirect product. Its normal factor is closed and
simply connected and its quotient is the simply connected abelian
group with Lie algebra $\mathfrak h_i$. Total geodesy follows from
the strong ideal property.
\end{proof}

\subsection[The isotropic right annihilator]{The case when \texorpdfstring{$Ann_R(\mathfrak g)$}{the right annihilator} is isotropic}
\label{subsect3}
Throughout this subsection, assume that $Ann_R(\mathfrak g)$ is
isotropic.

\begin{proposition}\label{prop4}
Under this assumption, if
$\mathfrak g=\mathfrak m_1\oplus\mathfrak m_2$ is a direct sum of
strong ideals, the restrictions of the metric to both summands
are non-degenerate.
\end{proposition}
\begin{proof}
Let $X\in\mathfrak m_1\cap\mathfrak m_1^\perp$.
For $Y\in\mathfrak m_1$ and $Z\in\mathfrak g$,
$$\langle\nabla_YX,Z\rangle=-\langle X,\nabla_YZ\rangle=0.$$
The derivatives in directions in $\mathfrak m_2$ vanish by
Proposition~\ref{prop:mixed}. Hence $X\in Ann_R(\mathfrak g)$.
By isotropy and Proposition~\ref{prop2.1},
$$X\in Ann_R(\mathfrak g)^\perp
=\nabla_{\mathfrak g}\mathfrak g
=\nabla_{\mathfrak m_1}\mathfrak m_1
\oplus\nabla_{\mathfrak m_2}\mathfrak m_2.$$
Since $X\in\mathfrak m_1$, it belongs to the first product space.
This space is orthogonal to $\mathfrak m_2$ by
Proposition~\ref{prop:mixed}. Thus $X$ is orthogonal to all of
$\mathfrak g$, and $X=0$. The other summand is treated in the
same way.
\end{proof}

Let
\begin{equation}\label{eq:two}
\mathfrak g=\bigoplus_{i=1}^n\mathfrak g_i
=\bigoplus_{j=1}^m\mathfrak g'_j
\end{equation}
be decompositions into indecomposable, non-degenerate strong
ideals. The right annihilator of each summand is contained in
$Ann_R(\mathfrak g)$ and is therefore isotropic. In particular,
$\nabla_{\mathfrak g_i}\mathfrak g_i\ne0$. After reordering the
second decomposition, we may assume
$\nabla_{\mathfrak g_1}\mathfrak g'_1\ne0$. Put
$$\mathfrak h_1=\bigoplus_{i=2}^n\mathfrak g_i,
\qquad \mathfrak h'_1=\bigoplus_{j=2}^m\mathfrak g'_j.$$

\begin{lemma}\label{lemma1}
$\mathfrak g_1\cap\mathfrak h'_1=0$ and
$\mathfrak g'_1\cap\mathfrak h_1=0$.
\end{lemma}
\begin{proof}
Write $A=\mathfrak g_1$, $B_1=A\cap\mathfrak g'_1$ and
$B_2=A\cap\mathfrak h'_1$. Then $B_1\ne0$. Let $p$ be the
projection onto $A$ in the first decomposition and $q$ the
projection onto $\mathfrak g'_1$ in the second. Both commute with
all left and right multiplications. Set
$$U=pq(A),\qquad V=p(I-q)(A).$$
We have $A=U+V$, $B_1\subseteq U$, $B_2\subseteq V$, and
\begin{align*}
\nabla_AU+\nabla_UA&\subseteq B_1,\\
\nabla_AV+\nabla_VA&\subseteq B_2.
\end{align*}
For example, if $x,y\in A$, then
$\nabla_y(pqx)=\nabla_y(qx)=q(\nabla_yx)\in B_1$;
the other inclusions follow in the same way.
\medskip

Choose a basis of $A/(B_1\oplus B_2)$ from the images of vectors
in $U$ or $V$, and lift these basis vectors to the respective
spaces. Adjoin the lifts from $U$ to $B_1$ and those from $V$ to
$B_2$, obtaining subspaces $C_1,C_2$. The inclusions just proved
show that both are strong ideals of $A$, and the choice of basis
gives $A=C_1\oplus C_2$. By Proposition~\ref{prop4}, applied to
$A$ with its isotropic right annihilator, both are non-degenerate.
Indecomposability and $B_1\ne0$ imply $C_2=0$, so $B_2=0$.
Interchanging the two decompositions proves the second assertion.
The same argument applies whenever the two selected non-zero
product factors have isotropic right annihilators, even if the
complements contain an additional zero-product summand.
\end{proof}

\begin{lemma}\label{lemma2}
The projection $\pi_1:\mathfrak g_1\longrightarrow\mathfrak g'_1$
along $\mathfrak h'_1$ is a strong isometry.
\end{lemma}
\begin{proof}
Lemma~\ref{lemma1} makes the projections in both directions
injective. The dimensions are therefore equal and $\pi_1$ is a
strong isomorphism. Write $X=X_1+X_2\in\mathfrak g_1$ with
$X_1\in\mathfrak g'_1$ and $X_2\in\mathfrak h'_1$. The lemma
and the mixed-product identities give $X_2\in Ann_R(\mathfrak g)$.
For example,
$$\nabla_{\mathfrak h'_1}X_2
=\nabla_{\mathfrak h'_1}X
\subseteq\mathfrak h'_1\cap\mathfrak g_1=0.$$
The derivatives from $\mathfrak g'_1$ vanish as well.
\medskip

By Proposition~\ref{prop4}, $\mathfrak h'_1$ is non-degenerate.
Write $X_1=Y_1+Y_2$ according to its orthogonal splitting with
$Y_1\in\mathfrak h'_1$ and $Y_2\in(\mathfrak h'_1)^\perp$.
Proposition~\ref{prop1} gives
$\nabla_{\mathfrak h'_1}Y_2=0$, whence
$Y_1\in Ann_R(\mathfrak h'_1)\subseteq Ann_R(\mathfrak g)$.
By isotropy,
$$\langle X_1,X_2\rangle=\langle Y_1,X_2\rangle=0,
\qquad\langle X_2,X_2\rangle=0.$$
Thus $\langle X,X\rangle=\langle\pi_1X,\pi_1X\rangle$.
Polarization proves the assertion.
\end{proof}

The intersection argument also gives vanishing of both orders of
mixed multiplication between unmatched factors. The difference
$X-\pi_1X$ in the preceding proof then belongs to the two-sided
annihilator. Hence matched factors have the same product space.
Repeating the argument matches every factor with exactly one factor
of the other decomposition: a second match would contradict the
vanishing just obtained. No non-zero factor can remain unmatched,
since its product space is non-zero. We have proved the following.

\begin{theorem}\label{th2}
Suppose $Ann_R(\mathfrak g)$ is isotropic. Let (\ref{eq:two})
be two decompositions into indecomposable, non-degenerate strong
ideals. Then $n=m$. After reordering,
\begin{align*}
\nabla_{\mathfrak g_i}\mathfrak g_i
&=\nabla_{\mathfrak g_i}\mathfrak g'_i
=\nabla_{\mathfrak g'_i}\mathfrak g_i
=\nabla_{\mathfrak g'_i}\mathfrak g'_i,\\
\nabla_{\mathfrak g_i}\mathfrak g'_j
&=\nabla_{\mathfrak g'_j}\mathfrak g_i=0\quad(i\ne j).
\end{align*}
The factor projections $\pi_i:\mathfrak g_i\longrightarrow
\mathfrak g'_i$ are strong isometries, and their direct sum is a
strong automorphism of $\mathfrak g$.
\end{theorem}

\begin{remark}
The last strong automorphism is not asserted to preserve the
metric between different factors. Factorwise isometry and a
global strong isometry are different conclusions when the
decompositions are not orthogonal.
\end{remark}

\subsection[Equality of the two annihilators]{The case when \texorpdfstring{$Ann_R(\mathfrak g)=Ann(\mathfrak g)$}{the right annihilator}}
\label{subsect4}
Suppose first that $\mathfrak g=\mathfrak g_1\oplus\mathfrak g_2$
is a non-degenerate strong decomposition. For
$X\in\mathfrak g_1^\perp$, write $X=X_1+X_2$ in this direct sum.
For $Y,Z\in\mathfrak g_1$,
$$\langle\nabla_YX_1,Z\rangle
=-\langle X_1,\nabla_YZ\rangle
=\langle X_2,\nabla_YZ\rangle=0.$$
It follows that $X_1\in Ann_R(\mathfrak g)=Ann(\mathfrak g)$.
Consequently $\nabla_X\mathfrak g_1=0$.
Together with Proposition~\ref{prop1}, this shows that
$\mathfrak g_1^\perp$ is a non-degenerate strong ideal. Thus a
decomposable algebra in this case also admits an orthogonal
decomposition.
\medskip

Choose a maximal non-degenerate subspace $\mathfrak g_0$ of
$Ann(\mathfrak g)$. Both $\mathfrak g_0$ and
$\mathfrak g_0^\perp$ are strong ideals. The right annihilator
of $\mathfrak g_0^\perp$ is isotropic: otherwise it would contain
a non-isotropic vector which could be added to $\mathfrak g_0$.
Splitting the complement orthogonally now gives the following.

\begin{proposition}\label{orth}
If $Ann_R(\mathfrak g)=Ann(\mathfrak g)$, there exists an
orthogonal decomposition
$$\mathfrak g=\mathfrak g_0\oplus\mathfrak g_1\oplus\cdots
\oplus\mathfrak g_l$$
where $\mathfrak g_0$ is a maximal non-degenerate subspace of
$Ann(\mathfrak g)$, and each $\mathfrak g_i$, $i\geq1$, is an
indecomposable non-degenerate strong ideal with isotropic right
annihilator. The summand $\mathfrak g_0$ is allowed to be zero.
\end{proposition}

\begin{exam}\label{rem2}
The assumption in Proposition~\ref{orth} cannot simply be omitted.
Let $\mathfrak g=\mathfrak e_1\oplus\mathfrak e_2$, with each
$\mathfrak e_i$ a copy of Example~\ref{ex:e2}, having basis
$\{a_i,u_i,v_i\}$. Give these six vectors unit length and add
only the mixed scalar product
$$\langle a_1,a_2\rangle=\rho,\qquad 0<|\rho|<1.$$
This metric is positive definite. Since the product space in
$\mathfrak e_i$ is $\operatorname{span}\{u_i,v_i\}$,
Proposition~\ref{prop:coupling} shows that the only non-zero
products are
$$\nabla_{a_i}u_i=v_i,\qquad\nabla_{a_i}v_i=-u_i
\quad(i=1,2).$$
Both summands are indecomposable strong ideals, and the metric is
flat. Moreover,
$$Ann_R(\mathfrak g)=\operatorname{span}\{a_1,a_2\},
\qquad Ann(\mathfrak g)=0.$$
\medskip

There is no alternative non-trivial orthogonal strong decomposition.
Indeed, any projection satisfying (\ref{eq:projector}) preserves
each plane $U_i=\operatorname{span}\{u_i,v_i\}$ and their common
complement $\operatorname{span}\{a_1,a_2\}$, since it commutes
with $L_{a_1}$ and $L_{a_2}$. Write
$Pa_i=b_{1i}a_1+b_{2i}a_2$. For $j\ne i$ and $u\in U_j$,
$$0=P\nabla_{a_i}u=\nabla_{Pa_i}u=b_{ji}L_{a_j}u,$$
so $b_{ji}=0$. For $u\in U_i$ the same identity yields
$P|_{U_i}=b_{ii}I$. Therefore
$$P=\alpha I_{\mathfrak e_1}\oplus\beta I_{\mathfrak e_2},
\qquad\alpha,\beta\in\{0,1\}.$$
The only non-trivial splitting is the displayed one, up to order,
and its mixed scalar product is $\rho\ne0$. Thus decomposability
does not imply orthogonal decomposability, even for a flat
Riemannian metric. At the same time this example shows that
$Ann_R(\mathfrak g)=0$ is sufficient, but not necessary, for
uniqueness of the actual strong factors.
\end{exam}

We next prove uniqueness for the orthogonal decompositions in
Proposition~\ref{orth}. The correction of the factor projections
is relevant when the non-degenerate annihilator summand is non-zero.

\begin{theorem}\label{th1}
Suppose $Ann_R(\mathfrak g)=Ann(\mathfrak g)$ and
$$\mathfrak g=\mathfrak g_0\oplus\bigoplus_{i=1}^n\mathfrak g_i
=\mathfrak g'_0\oplus\bigoplus_{j=1}^m\mathfrak g'_j$$
are orthogonal decompositions as in Proposition~\ref{orth}.
Then $n=m$. After reordering the non-zero-product factors, their
product spaces agree and their mixed products with unmatched
factors vanish as in Theorem~\ref{th2}. There exist strong
isometries $\phi_i:\mathfrak g_i\longrightarrow\mathfrak g'_i$,
$0\leq i\leq n$, whose direct sum is a strong isometry of
$\mathfrak g$. For $i\geq1$, $\phi_i$ can be chosen to be the
identity on $D_i=\nabla_{\mathfrak g_i}\mathfrak g_i$.
\end{theorem}
\begin{proof}
The intersection argument in Lemma~\ref{lemma1} applies to the
factors with non-zero product space, since their right annihilators
are isotropic. It matches these factors bijectively, gives
$D_i=D'_i$, and makes the restricted projections
$q_i:\mathfrak g_i\longrightarrow\mathfrak g'_i$ strong
isomorphisms fixing $D_i$. It also gives vanishing of the products
between unmatched factors. Write
$A_i=Ann_R(\mathfrak g_i)=Ann(\mathfrak g_i)$. Since $A_i$ is
isotropic, Proposition~\ref{prop2.1} gives $A_i\subseteq D_i$.
As $A_i\subseteq Ann(\mathfrak g)$ and $D_i=D'_i$, we obtain
$A_i=A'_i$.
\medskip

Let $q_0$ be the projection onto $\mathfrak g'_0$, restricted to
$\mathfrak g_i$. The projection of $\mathfrak g_i$ into any
unmatched non-zero-product factor is contained in its isotropic
annihilator. Orthogonality therefore gives, for $X,Y\in\mathfrak g_i$,
\begin{equation}\label{eq:defect}
\langle q_iX,q_iY\rangle
=\langle X,Y\rangle-\delta(X,Y),
\qquad\delta(X,Y)=\langle q_0X,q_0Y\rangle.
\end{equation}
The symmetric form $\delta$ vanishes whenever either argument
belongs to $D_i$. Choose a basis $a_1,\ldots,a_k$ of $A_i$,
a non-degenerate complement $E_i$ of $A_i$ in $D_i$, and vectors
$b_1,\ldots,b_k$ so that
$$\langle a_r,b_s\rangle=\delta_{rs},\quad
\langle b_r,b_s\rangle=0,\quad
E_i\perp\operatorname{span}\{a_r,b_r:1\leq r\leq k\}.$$
Define $\phi_i$ to fix $D_i$ and set
$$\phi_i(b_r)=q_i(b_r)
+\frac12\sum_{s=1}^k\delta(b_r,b_s)a_s.$$
The added vectors belong to $Ann(\mathfrak g)$, and the correction
vanishes on the product space. Hence $\phi_i$ remains a strong
isomorphism. Formula (\ref{eq:defect}) shows that it preserves all
pairings of these basis vectors, so it is an isometry.
\medskip

The two maximal non-degenerate annihilator summands have the same
signature, namely the signature of the non-degenerate part of
the restricted form on $Ann(\mathfrak g)$. Choose any linear
isometry $\phi_0$ between them; it is strong because their products
are zero. Orthogonality of both decompositions makes
$\bigoplus_{i=0}^n\phi_i$ a strong isometry of $\mathfrak g$.
\end{proof}

\begin{remark}\label{rem:projectioncorrection}
The maps $\phi_i$ in Theorem~\ref{th1} need not be the uncorrected
projections $q_i$. Their metric defect is exactly
(\ref{eq:defect}). If the non-degenerate annihilator part vanishes,
the correction is unnecessary. In general it cannot be omitted.
\end{remark}

\subsection[The zero right annihilator]{The case when \texorpdfstring{$Ann_R(\mathfrak g)=0$}{the right annihilator}}\label{subsect5}
\begin{proposition}\label{prop6}
If $Ann_R(\mathfrak g)=0$, every direct sum decomposition into
strong ideals is non-degenerate and orthogonal. The following
conditions are equivalent:
\begin{enumerate}
\item $\mathfrak g$ is decomposable;
\item $\mathfrak g$ is a direct sum of two non-zero strong ideals;
\item $\mathfrak g$ is an orthogonal direct sum of two non-zero
Lie ideals.
\end{enumerate}
\end{proposition}
\begin{proof}
By Proposition~\ref{prop2.1}, $\mathfrak g=\nabla_{\mathfrak g}
\mathfrak g$. For a strong direct sum, mixed products vanish and
each factor satisfies $\mathfrak g_i=\nabla_{\mathfrak g_i}
\mathfrak g_i$. Proposition~\ref{prop:mixed} then makes distinct
factors orthogonal; they are non-degenerate as well. Conversely,
the Koszul formula shows that an orthogonal direct sum of Lie
ideals is a strong direct sum. This last implication does not
require the annihilator assumption.
\end{proof}

\begin{theorem}\label{th3}
If $Ann_R(\mathfrak g)=0$, there is an orthogonal decomposition
$$\mathfrak g=\mathfrak g_1\oplus\cdots\oplus\mathfrak g_l$$
into indecomposable, non-degenerate strong ideals, each with zero
right annihilator. The factors, as subspaces of $\mathfrak g$,
are unique up to their order.
\end{theorem}
\begin{proof}
Existence and orthogonality follow from Propositions~\ref{prop2}
and \ref{prop6}. We give a direct proof of uniqueness. If
$\mathfrak g=\bigoplus_j\mathfrak g'_j$ is another such
decomposition, then
$$\mathfrak g_i=\nabla_{\mathfrak g_i}\mathfrak g
=\bigoplus_j\nabla_{\mathfrak g_i}\mathfrak g'_j
\subseteq\bigoplus_j(\mathfrak g_i\cap\mathfrak g'_j)
\subseteq\mathfrak g_i.$$
Consequently $\mathfrak g_i$ is the direct sum of these
intersections. They are strong ideals and mutually orthogonal,
since the second decomposition is orthogonal. Non-degeneracy of
$\mathfrak g_i$ implies non-degeneracy of each non-zero
intersection. Indecomposability leaves only one, so
$\mathfrak g_i\subseteq\mathfrak g'_j$ for some $j$.
Applying the same argument to $\mathfrak g'_j$ gives equality.
All factors are therefore paired, proving uniqueness.
\end{proof}

\section{Some applications of the decomposition results}\label{sect4}
We use the curvature convention
$$R(X,Y)Z=\nabla_X\nabla_YZ-\nabla_Y\nabla_XZ-\nabla_{[X,Y]}Z,$$
and define the Ricci tensor by
$$ric(X,Y)=\operatorname{tr}\{Z\longmapsto R(Z,X)Y\}.$$
The metric is {\bf flat} if $R=0$, {\bf Einstein} if
$ric=c\langle,\rangle$ for a constant $c$, and {\bf Ricci flat}
if $ric=0$.

\subsection{Einstein metrics and canonical factors}
\begin{proposition}\label{prop5}
$Ann_R(\mathfrak g)\subseteq\ker ric$. In particular, if $ric$ is
non-degenerate, then $Ann_R(\mathfrak g)=0$ and
$\mathfrak g=\nabla_{\mathfrak g}\mathfrak g$.
\end{proposition}
\begin{proof}
If $X\in Ann_R(\mathfrak g)$, each term in $R(Z,Y)X$ vanishes.
Taking the trace gives $ric(Y,X)=0$ for every $Y$. The remaining
assertions follow from non-degeneracy and Proposition~\ref{prop2.1}.
\end{proof}

\begin{theorem}\label{th4}
Let $G$ have a left invariant pseudo-Riemannian metric with
non-degenerate Ricci tensor. There is a decomposition
$$\mathfrak g=\mathfrak g_1\oplus\cdots\oplus\mathfrak g_l$$
into indecomposable, non-degenerate strong ideals, unique up to
their order. Distinct factors are orthogonal for both the metric
and the Ricci tensor. The Ricci tensor restricted to a factor is
the Ricci tensor of the restricted metric, and is non-degenerate.
\end{theorem}
\begin{proof}
Proposition~\ref{prop5} and Theorem~\ref{th3} give the strong
decomposition, its orthogonality and its uniqueness. By
Proposition~\ref{prop:mixed}, if $X=\sum_iX_i$, $Y=\sum_iY_i$
and $Z=\sum_iZ_i$, then
$$R(X,Y)Z=\sum_iR^i(X_i,Y_i)Z_i.$$
Taking the trace gives
$$ric(X,Y)=\sum_i ric_i(X_i,Y_i).$$
Thus all mixed Ricci terms vanish and the restrictions have the
asserted form. Their non-degeneracy follows from that of $ric$.
\end{proof}

\begin{corollary}\label{cor1}
If $\langle,\rangle$ is Einstein with constant $c\ne0$, its
indecomposable strong factors are orthogonal and unique up to
their order. Each restricted metric is Einstein with the same
constant $c$.
\end{corollary}

\begin{proposition}\label{prop:einsteinconstruction}
Let $(\mathfrak g_i,\langle,\rangle_i)$ be Einstein metric Lie
algebras with constants $c_i$, and let $a_i>0$. The orthogonal
direct sum metric $\bigoplus_i a_i\langle,\rangle_i$ is Einstein
if and only if all numbers $c_i/a_i$ agree.
\end{proposition}
\begin{proof}
Multiplication of a metric by a non-zero constant does not change
its Levi-Civita connection or its Ricci tensor as a covariant
two-tensor. On the $i$th factor we therefore have
$$ric_i=c_i\langle,\rangle_i
=\frac{c_i}{a_i}\big(a_i\langle,\rangle_i\big).$$
The Ricci tensor of an orthogonal product is the direct sum of
these tensors. This proves both directions.
\end{proof}

\begin{exam}\label{ex:einstein}
Let $\mathfrak a$ have orthonormal basis $\{a,b\}$ and bracket
$[a,b]=b$. Its only non-zero products are
$$\nabla_b a=-b,\qquad\nabla_b b=a,$$
so $R(a,b)b=-a$ and $ric=-\langle,\rangle$. Two copies of this
metric with the same positive scale give an Einstein metric on
the product with constant $-1/a_1$. Unequal scales give a metric
with non-degenerate Ricci tensor which is not Einstein.
The two affine Lie algebra factors are strongly indecomposable,
since their centers are zero and their dimensions are two.
Their product decomposition is therefore the unique one in
Theorem~\ref{th4}.
\end{exam}

\begin{corollary}\label{cor:automorphisms}
Assume $Ann_R(\mathfrak g)=0$. Every isometric Lie algebra
automorphism permutes the indecomposable strong factors. The
identity component of the group of such automorphisms preserves
each factor. Consequently every skew-adjoint derivation preserves
each factor. These conclusions apply in particular when $ric$ is
non-degenerate.
\end{corollary}
\begin{proof}
An isometric Lie algebra automorphism preserves the Levi-Civita
product by (\ref{nabla1}), so it sends a strong decomposition to
another one. Theorem~\ref{th3} gives a permutation of the factors.
A continuous path starting at the identity induces a constant
permutation of this finite set. Exponentiating a skew-adjoint
derivation and differentiating gives the last assertion.
\end{proof}

This gives a reduction of the symmetry calculation as well as of
the Einstein equation. It concerns automorphisms compatible with
the specified Lie group structure. It makes no assertion that an
arbitrary isometry of the underlying manifold is a group
automorphism. Likewise, Example~\ref{ex:compact} shows why an
arbitrary algebraic direct sum cannot be substituted for a strong
decomposition in these applications.

\subsection{The Ricci-flat boundary}
\begin{exam}\label{exam1}
Suppose the metric is bi-invariant and let $K$ be the Killing form.
Then (\ref{nabla1}) gives
$$\nabla_XY=\tfrac12[X,Y],\qquad
R(X,Y)Z=-\tfrac14[[X,Y],Z],\qquad
ric(X,Y)=-\tfrac14K(X,Y).$$
In particular the metric is flat if and only if
$[\mathfrak g,[\mathfrak g,\mathfrak g]]=0$, that is, the algebra
is nilpotent of step at most two. If it is Einstein with constant
$c\ne0$, the Killing form is non-degenerate and the algebra is
semisimple. If it is Ricci flat, then $K=0$ and Cartan's solvability
criterion makes the algebra solvable. A nilpotent algebra with an
invariant non-degenerate form and nilpotency step at least three
therefore gives a non-flat Ricci-flat example.
\end{exam}

\begin{exam}\label{ex:nilpotent}
For a concrete instance, let $\mathfrak n$ have basis
$e_1,\ldots,e_5$ and brackets
$$[e_1,e_2]=e_3,\qquad [e_1,e_3]=e_4,\qquad [e_2,e_3]=e_5.$$
Give it the scalar products
$$\langle e_1,e_5\rangle=1,\qquad
\langle e_2,e_4\rangle=-1,\qquad\langle e_3,e_3\rangle=1.$$
The bracket satisfies the Jacobi identity and defines a three-step
nilpotent Lie algebra. The form is non-degenerate of signature
$(3,2)$. Pairing the displayed brackets with the basis verifies
$\langle[X,Y],Z\rangle+\langle Y,[X,Z]\rangle=0$, so it is
invariant. Its Killing form is zero, but
$$R(e_1,e_2)e_1=\tfrac14e_4\ne0.$$
Thus the metric on the associated simply connected group is
bi-invariant and Ricci flat, but is not flat. Here
$$Ann_R(\mathfrak n)=Ann(\mathfrak n)
=\operatorname{span}\{e_4,e_5\}$$
is totally isotropic.
\end{exam}

In contrast, a homogeneous Ricci-flat Riemannian metric is flat
\cite{AK}. Flat left invariant pseudo-Riemannian metrics are
discussed, for example, in \cite{Ba,CHB}. For the Riemannian case
we recall the following precise structure theorem.

\begin{theorem}[\cite{Mi}, Theorem 1.5]\label{th:milnor}
A Lie group with a left invariant Riemannian metric is flat if and
only if its metric Lie algebra is an orthogonal vector space sum
$\mathfrak b\oplus\mathfrak u$, where $\mathfrak b$ is an abelian
subalgebra, $\mathfrak u$ is an abelian ideal, and
$\operatorname{ad}b$ is skew-adjoint for every $b\in\mathfrak b$.
In this case $\nabla_u=0$ and $\nabla_b=\operatorname{ad}b$.
\end{theorem}

\begin{theorem}\label{th6}
Equivalently, after separating the common fixed subspace of the
skew-adjoint action in Theorem~\ref{th:milnor}, a flat Riemannian
metric Lie algebra has an orthogonal vector space sum
$$\mathfrak g=\mathfrak b\oplus Ann(\mathfrak g)
\oplus[\mathfrak g,\mathfrak g].$$
Here $\mathfrak b$ is an abelian subalgebra,
$[\mathfrak g,\mathfrak g]=\nabla_{\mathfrak g}\mathfrak g$ is
an abelian ideal of even dimension, and
\begin{align*}
Ann_R(\mathfrak g)&=\mathfrak b\oplus Ann(\mathfrak g),\\
\dim\mathfrak b&\leq\tfrac12\dim[\mathfrak g,\mathfrak g].
\end{align*}
For $b\in\mathfrak b$ we have $\nabla_Xb=0$ for every $X$, and
$\nabla_b=\operatorname{ad}b$. The latter action is injective on
$\mathfrak b$ and skew-adjoint.
\end{theorem}
\begin{proof}
The common fixed vectors in $\mathfrak u$, together with any
kernel of the action of the original $\mathfrak b$, are central
and belong to $Ann(\mathfrak g)$. Separate them orthogonally.
The remaining commuting skew-adjoint operators act faithfully
on the orthogonal complement of their common kernel. Simultaneous
orthogonal block decomposition gives two-dimensional rotation
planes. Their images span the remaining abelian ideal, which is
both $[\mathfrak g,\mathfrak g]$ and
$\nabla_{\mathfrak g}\mathfrak g$. Its dimension is even and
faithfulness bounds $\dim\mathfrak b$ by the number of rotation
planes. The annihilator formulas follow from the displayed
connection. This is also the description in \cite{CHB}.
\end{proof}

\begin{remark}\label{rem:flatunique}
For a non-zero flat Riemannian metric Lie algebra the right
annihilator is non-zero. Uniqueness of the factors as actual
subspaces is not guaranteed: in the abelian Euclidean plane every
orthogonal pair of lines is a strong decomposition, and rotations
give distinct pairs. They are nevertheless related by strong
isometries. This shows that the conclusion ``unique up to order''
in Corollary~\ref{cor1} cannot be extended to all Ricci-flat
metrics. It does not mean that every Ricci-flat example lacks
uniqueness; Example~\ref{rem2} has unique strong factors.
\end{remark}

\subsection{Geodesics and completeness}
\begin{proposition}\label{prop:geodesics}
Let $\mathfrak g=\bigoplus_{i=1}^l\mathfrak g_i$ be a
non-degenerate strong decomposition and let
$\widetilde G=\prod_i\widetilde G_i$ be its simply connected
group decomposition. A curve
$\gamma(t)=(\gamma_1(t),\ldots,\gamma_l(t))$ is an affinely
parametrized geodesic if and only if every $\gamma_i(t)$ is an
affinely parametrized geodesic of the restricted metric, with the
same parameter. Moreover $\widetilde G$ is geodesically complete
if and only if every $\widetilde G_i$ is geodesically complete.
The decomposition need not be orthogonal.
\end{proposition}
\begin{proof}
By Proposition~\ref{prop:mixed}, the connection under the product
identification is the direct sum of the restricted connections.
More explicitly, the left translated velocity
$$u(t)=(L_{\gamma(t)^{-1}})_*\dot\gamma(t)
=\sum_i u_i(t)$$
satisfies the geodesic equation
$$\dot u+\nabla_u u=0.$$
The vanishing of mixed products makes this equivalent to
\begin{equation}\label{eq:geodesicparts}
\dot u_i+\nabla^i_{u_i}u_i=0,\qquad 1\leq i\leq l.
\end{equation}
Thus complete factors give a geodesic through every initial
condition for all real times. Conversely, take zero initial
velocities in all but one factor. A global geodesic in the product
then supplies a global geodesic in the remaining factor.
\end{proof}

\begin{corollary}\label{cor:completecover}
For any connected $G$ with the lifted strong decomposition on
$\widetilde G$, the metric on $G$ is geodesically complete if and
only if each simply connected factor $\widetilde G_i$ is
geodesically complete.
\end{corollary}
\begin{proof}
A pseudo-Riemannian covering is a local isometry. Geodesics lift
and project with their affine parameters, so completeness is
equivalent on the base and its covering group. Apply
Proposition~\ref{prop:geodesics}.
\end{proof}

\begin{exam}\label{ex:geodesicflat}
For the non-orthogonal product in Example~\ref{rem2}, write
$$\xi_i(t)=c_i(t)a_i+p_i(t)u_i+q_i(t)v_i.$$
Equation (\ref{eq:geodesicparts}) becomes
$$\dot c_i=0,\qquad\dot p_i=c_iq_i,
\qquad\dot q_i=-c_ip_i.$$
Hence $c_i$ is constant and $p_i+\sqrt{-1}q_i$ rotates with
constant angular speed $-c_i$. These equations are independent
of the mixed metric parameter $\rho$. All factors and the product
are complete. Their total energy, however, includes the mixed term:
$$\langle u(t),u(t)\rangle
=\sum_{i=1}^2(c_i^2+p_i^2+q_i^2)+2\rho c_1c_2.$$
This is a concrete use of the non-orthogonal decomposition: the
geodesic equations separate even though the metric energy is not
the sum of the restricted energies alone.
\end{exam}

\begin{remark}
If $V\in Ann_R(\mathfrak g)$, every geodesic has a conserved
quantity $\langle\dot\gamma,V\rangle$, since $V$ is parallel.
Under a strong product, the factor energies are also conserved.
Proposition~\ref{prop:coupling} explains why allowed mixed metric
terms are constant along solutions of the separated equations.
\end{remark}

\section{Further remarks}\label{sect:meaning}
The preceding results distinguish three kinds of information.
An ordinary Lie algebra decomposition concerns the bracket; the
de Rham--Wu decomposition concerns parallel non-degenerate
distributions of the metric manifold; a strong decomposition
requires the Levi-Civita product to split into two-sided ideals.
It therefore gives normal group factors whose mixed covariant
derivatives vanish. An orthogonal strong decomposition is a metric
product compatible with group multiplication. A non-orthogonal
one is still a product for the connection, with its remaining
metric terms described by Proposition~\ref{prop:coupling}.
\medskip

Examples~\ref{ex:e2} and \ref{ex:twisted} show that metric product
structure alone does not recover the strong factors, even in
positive definite signature and even when $Ann_R(\mathfrak g)=0$.
In particular $Ann_R(\mathfrak g)$ measures left invariant
parallel fields, not all parallel fields of the manifold. It must
not be identified with the fixed subspace of the full holonomy
representation without an additional argument about left
invariance.
\medskip

The product spaces $D_i=\nabla_{\mathfrak g_i}\mathfrak g_i$
give a useful explanation of the uniqueness results. Mixed metric
terms vanish on these spaces, and matched factors in the isotropic
case have the same $D_i$. If $Ann_R(\mathfrak g)=0$, then
$D_i=\mathfrak g_i$, so neither mixed terms nor changes of the
factor subspaces survive. If annihilators remain, the freedom of
mixing factors is constrained by the quotient spaces
$\mathfrak g_i/D_i$ and the annihilator pairings. Theorems
\ref{th2} and \ref{th1} describe two forms of uniqueness in this
setting; they do not give an unrestricted uniqueness theorem for
every possible annihilator.
\medskip

For Einstein metrics with non-zero constant, the decomposition
reduces the equation and the compatible automorphism group to
indecomposable factors. For arbitrary signature, the geodesic
application reduces a global completeness question to the
individual factors without assuming a metric product. These
consequences use exactly the extra compatibility with the
Levi-Civita product.
\medskip

Two further questions are suggested by the examples. One is to
find sharper conditions for uniqueness when
$Ann_R(\mathfrak g)\ne0$, since this condition alone does not
decide uniqueness. Another is to describe strongly indecomposable
metric Lie groups whose underlying metric manifolds admit a
non-trivial de Rham--Wu product, beyond the flat and non-flat
examples above. These are questions for further study, rather
than consequences of the present uniqueness theorems.
\section{Acknowledgments}
This work is supported by National Natural Science Foundation of China (No. 11001133) and the Fundamental Research Funds for the Central Universities. This work was partly completed when the first author visited the University of California, Berkeley. The first author would like to thank Joseph A. Wolf for the invitation, helpful discussions and suggestions.


\begin{thebibliography}{99}
\bibitem{AK}
\textsc{D. Alekseevsky} and \textsc{B. Kimel'fel'd}, Structure of homogeneous Riemannian spaces with zero Ricci curvature, \textit{Functional Anal. Appl.}
\textbf{9} (1975), 97--102.

\bibitem{An}
\textsc{M.T. Anderson}, Einstein metrics and metrics with bounds on Ricci curvature, \textit{Proc. I.C.M.} (1994, Zurich), 443--452, Basel, Birkh\"auser, 1995.

\bibitem{Ba}
\textsc{O. Baues}, Prehomogeneous affine representations and flat pseudo-Riemannian manifolds, Handbook of pseudo-Riemannian geometry and supersymmetry, 731--817, IRMA Lect. Math. Theor. Phys., 16, Eur. Math. Soc., Zurich, 2010.

\bibitem{Ber}
\textsc{M. Berger}, Sur Les groupes d'holomomie des varietes a connexion affine et des varietes riemanniennes, \textit{Bull. Soc. Math. France} \textbf{83} (1955), 279--330.

\bibitem{Ber1}
\textsc{M. Berger}, Riemannian geometry during the second half of the twentieth century, \textit{Univ. Lect. Ser.} \textbf{17}, Amer. Math. Soc., Providence, 1998.

\bibitem{Ber2}
\textsc{M. Berger}, A panoramic view of Riemannian geometry, Springer-Verlag, Berlin-Heidelberg, 2003.

\bibitem{Be}
\textsc{A. Besse}, Einstein manifolds, \textit{Ergeb. Math.} \textbf{10} (1987), Springer-Verlag, Berlin-Heidelberg.

\bibitem{Bo}
\textsc{C. B$\rm{\ddot{o}}$hm}, Homogeneous Einstein metrics and simplicial complexes, \textit{J. Differential Geom}. \textbf{67} (2004), no. 1, 79--165.

\bibitem{Ch}
\textsc{Z.Q. Chen}, The uniqueness in the de Rham--Wu decomposition, arXiv:1104.4165 [math.DG]. \url{https://arxiv.org/abs/1104.4165}.

\bibitem{CHB}
\textsc{Z.Q. Chen}, \textsc{D.P. Hou} and \textsc{C.M. Bai}, A
left-symmetric algebraic approach to left invariant flat
(pseudo-)metrics on Lie groups, \textit{J. Geom. Phys.}
\textbf{62} (2012), no. 7, 1600--1610.
\url{https://doi.org/10.1016/j.geomphys.2012.03.003}.

\bibitem{CLZ}
\textsc{Z.Q. Chen}, \textsc{K. Liang} and \textsc{F.H. Zhu}, Algebras with pseudo-Riemannian bilinear forms, \textit{J. Korean. Math. Soc.} \textbf{48} (2011), no. 1, 1--12.

\bibitem{DZ1}
\textsc{J.E. D'Atri} and \textsc{W. Ziller}, Naturally reductive metrics and Einstein metrics on compact Lie groups, \textit{Mem. Amer. Math. Soc.} \textbf{18} (1979), no. 215.

\bibitem{de}
\textsc{G. de Rham}, Sur la reductibilite d'un espace de Riemann, \textit{Comm. Math. Helv.} \textbf{26} (1952), 329--344.

\bibitem{GT}
\textsc{A.S. Galaev} and \textsc{T. Leistner}, Recent developments in pseudo-Riemannian holonomy theory. \textit{Handbook of pseudo-Riemannian geometry and supersymmetry}, 581--627, IRMA Lect. Math. Theor. Phys., 16, Eur. Math. Soc., Zurich, 2010.

\bibitem{GK}
\textsc{C.S. Gordon} and \textsc{M. Kerr}, New homogeneous Einstein metrics of negative Ricci curvature, \textit{Ann. Global Anal. Geom.} \textbf{19} (2001), no. 1, 75--101.

\bibitem{He}
\textsc{J. Heber}, Noncompact homogeneous Einstein spaces, \textit{Invent. Math.} \textbf{133} (1998), no. 2, 279--352.

\bibitem{Ja1}
\textsc{M. Jablonski}, Detecting orbits along subvarieties via the moment map, \textit{M$\ddot{u}$nster J. Math.} \textbf{3} (2010), 67--88.

\bibitem{Ja2}
\textsc{M. Jablonski}, Concerning the existence of Einstein and Ricci soliton metrics on solvable Lie groups, \textit{Geom. Topol.} \textbf{15} (2011), 735--764.

\bibitem{Kn1}
\textsc{A. Knapp}, Lie groups beyond an introduction, second edition, Birkh\"auser, 2002.

\bibitem{La1}
\textsc{J. Lauret}, Einstein solvmanifolds and nilsolitons, New developments in Lie theory and geometry, 1--35, \textit{Contemp. Math.} \textbf{491}, Amer. Math. Soc., Providence, RI, 2009.

\bibitem{La2}
\textsc{J. Lauret}, Einstein solvmanifolds are standard, \textit{Ann. of Math.} \textbf{(2)172} (2010), 1859--1877.

\bibitem{LW}
\textsc{C. Lebrun} and \textsc{M.Y. Wang}, Surveys in differential geometry: essays on Einstein manifolds, International Press, Boston, 1999.

\bibitem{Ma}
\textsc{R. Maltz}, The de Rham product decomposition, \textit{J. Differential Geom.} \textbf{7} (1972), 161--174.

\bibitem{Mi}
\textsc{J. Milnor}, Curvatures of left invariant metrics on Lie groups, \textit{Adv. Math.} \textbf{21} (1976), 293--329. \url{https://doi.org/10.1016/S0001-8708(76)80002-3}.

\bibitem{Ni1}
\textsc{Y. Nikolayevsky}, Einstein solvmanifolds and the pre-Einstein derivation, \textit{Trans. Amer. Math. Soc.} \textbf{363} (2011), no. 8, 3935--3958.

\bibitem{Ni2}
\textsc{Y. Nikolayevsky}, Einstein solvmanifolds attached to two-step nilradicals, arXiv:0805.0646 [math.DG]. \url{https://arxiv.org/abs/0805.0646}.

\bibitem{Wu1}
\textsc{H. Wu}, On the de Rham decomposition theorem, \textit{Ill. J. Math.} \textbf{8} (1964), 291--311. \url{https://doi.org/10.1215/ijm/1256059674}.

\bibitem{Wu2}
\textsc{H. Wu}, Holonomy groups of indefinite metrics, \textit{Pacific J. Math.} \textbf{20} (1967), 351--392.

\bibitem{ZZ}
\textsc{F.H. Zhu} and \textsc{L.S. Zhu}, The uniqueness of the decomposition quadratic Lie algebras, \textit{Comm. Algebra} \textbf{29} (2001), no. 11, 5145--5154.
\end{thebibliography}
\end{document}